\documentclass[12pt]{article}
\usepackage{amssymb}
\usepackage{mathrsfs}
\usepackage{cite}
\usepackage{tikz}
\usepackage{amsmath}
\usepackage{bm}
\usepackage{graphicx}
\usepackage{pst-poly}  
\usepackage{pst-plot}  
\usepackage{multirow}
\usepackage{bigstrut}

\newtheorem{thm}{Theorem}[section]
\newtheorem{lem}[thm]{Lemma}

\newtheorem{cor}[thm]{Corollary}

\newenvironment{pf}[1][Proof]{\noindent\textbf{#1.} }{\hfill\rule{1mm}{2mm}}

\makeatletter \@addtoreset{equation}{section} \makeatother

\begin{document}
\title{\bf Embedded connectivity of recursive networks\thanks{The work was supported by NNSF of China (10711233, 11571044),
 NSF of Hubei(2014CFB248), Young talent fund from Hubei EDU(Q20151311) and ORFIAM, Yangtze University(KF1505, KF1506).}}

\author
{Xiang-Jun Li$^a$ \quad Qi-Qi Dong$^a$\quad Zheng Yan$^a$ \quad
Jun-Ming Xu$^b$ \footnote{Corresponding author: xujm@ustc.edu.cn
(J.-M. Xu) }\\
\\
 {\small $^a$School of Information and Mathematics,}\\
 {\small Yangtze University, Jingzhou, Hubei, 434023, China}\\
 {\small $^b$School of Mathematical Sciences}\\
 {\small University of Science and Technology of China, Hefei, 230026, China}}
\date{}
\maketitle

\begin{abstract}

Let $G_n$ be an $n$-dimensional recursive network. The $h$-embedded
connectivity $\zeta_h(G_n)$ (resp. edge-connectivity $\eta_h(G_n)$)
of $G_n$ is the minimum number of vertices (resp. edges) whose
removal results in disconnected and each vertex is contained in an
$h$-dimensional subnetwork $G_h$. This paper determines $\zeta_h$
and $\eta_h$ for the hypercube $Q_n$ and the star graph $S_n$, and
$\eta_3$ for the bubble-sort network $B_n$.

\vskip6pt


\vskip6pt


\end{abstract}

\section{Introduction}

It is well known that interconnection networks play an important
role in parallel computing/communication systems. An interconnection
network can be modeled by a graph $G=(V, E)$, where $V$ is the set
of processors and $E$ is the set of communication links in the
network.

The {\it connectivity} $\kappa(G)$ (resp. {\it edge-connectivity}
$\lambda(G)$ ) of $G$ is defined as the minimum number of vertices
(resp. edges) whose removal from $G$ results in a disconnected
graph. The connectivity $\kappa(G)$ and edge-connectivity
$\lambda(G)$ of a graph $G$ are two important measurements for fault
tolerance of the network since the larger $\kappa(G)$ or
$\lambda(G)$ is, the more reliable the network is.

However, the definitions of $\kappa(G)$ and $\lambda(G)$ are
implicitly assumed that any subset of system components is equally
likely to be faulty simultaneously, which may not be true in real
applications, thus they underestimate the reliability of the
network. To overcome such a shortcoming, Harary~\cite{h83}
introduced the concept of conditional connectivity by appending some
requirements on connected components, Latifi {\it et
al.}~\cite{lhm94} specified requirements and proposed the concept of
the restricted $h$-connectivity. These parameters can measure fault
tolerance of an interconnection network more accurately than the
classical connectivity. The concepts stated here are slightly
different from theirs.

For a graph $G$, $\delta(G)$ denotes its minimum degree. A subset
$S\subset V(G)$ (resp. $F\subset E(G)$) is called an {\it
$h$-vertex-cut} (resp. {\it edge-cut}), if $G-S$ (resp. $G-F$) is
disconnected and $\delta(G-S)\geq h$. The {\it $h$-super
connectivity} $\kappa^{h}(G)$ (resp. {\it $h$-super
edge-connectivity} $\lambda^{h}(G)$) of $G$ is defined as the
cardinality of a minimum $h$-vertex-cut (resp. $h$-edge-cut) of $G$.

For any graph $G$ and any integer $h$, determining $\kappa^{h}(G)$
and $\lambda^{h}(G)$ is quite difficult, no polynomial algorithm to
compute them has been yet known so far. In fact, the existence of
$\kappa^{h}(G)$ and $\lambda^{h}(G)$ is an open problem for $h\geq
1$. Only a little knowledge of results have been known on
$\kappa^{h}$ and $\lambda^{h}$ for some special classes of graphs
for any $h$, such as the hypercube $Q_n$ and the star graph $S_n$.

In order to facilitate the expansion of the network, and to use the
same routing algorithm or maintenance strategy as used in the
original one, large-scale parallel computing systems always take
some networks of recursive structures as underlying topologies, such
as the hypercube $Q_n$, the star graph $S_n$, the bubble-sort graph
$B_n$ and so on. The presence of vertex and/or edge failures maybe
disconnects the entire network, one hopes that every remaining
component has undamaged subnetworks (i.e., smaller networks with
same topological properties as the original one). Under this
consideration, Yang {\it et al.}~\cite{yw12} proposed the concept of
embedded connectivity.

Let $G_n$ be an $n$-dimensional recursive network. For a positive
integer $h$ with $h\leq n-1$, there is a sub-network $G_h\subset
G_n$. Let $\delta_h=\delta(G_h)$.

A subset $F\subset V(G_n)$ (resp. $F\subset E(G_n)$) is an {\it
$h$-embedded vertex-cut} (resp. {\it $h$-embedded edge-cut}) if
$G_n-F$ is disconnected and each vertex is contained in an
$h$-dimensional subnetwork $G_h$. The {\it $h$-embedded
connectivity} $\zeta_h(G_n)$ (resp. {\it edge-connectivity}
$\eta_h(G_n)$) of $G_n$ is defined as the cardinality of a minimum
$h$-embedded vertex-cut (resp. $h$-embedded edge-cut) of $G_n$.

By definition, if $S$ is an $h$-embedded vertex-cut of $G_n$ with
$|S|=\zeta_h(G_n)$, then $G_n-S$ contains a sub-network $G_h$, and
so $\delta(G_n-S)\geq \delta_h$, which implies $S$ is a
$\delta_h$-vertex-cut of $G_n$. Thus, $\kappa^{\delta_h}(G_n)\leq
|S|=\zeta_h(G_n)$. Similarly, $\lambda^{\delta_h}(G_n)\leq
\eta_h(G_n)$. These facts are useful and we write them as the
following lemma.

 \begin{lem}\label{lem1.1}  For $h\leq n-1$,
$\zeta_h(G_n)\geq \kappa^{\delta_h}(G_n)$ if
$\kappa^{\delta_h}(G_n)$ exists, and
$\eta_h(G_n)\geq\lambda^{\delta_h}(G_n)$ if
$\lambda^{\delta_h}(G_n)$ exists.
\end{lem}

Using Lemma~\ref{lem1.1}, for a star graph $S_n$ and a bubble-sort
graph $B_n$, Yang {\it et al.}~\cite{yw12,yw14} determined
$\zeta_2(S_n)=2n-4$ for $n\ge 3$, $\eta_2(S_n)=2n-4$ for $n\ge 3$
and $\eta_3(S_n)=6(n-3)$ for $n\ge 4$; and $\zeta_2(B_n)=2n-4$ for
$n\ge 3$. In this paper, we will determine $\zeta_h$ and $\eta_h$
for $Q_n$ and $S_n$ for any $h\leq n-1$ and determine $\eta_3(B_n)$.

The rest of the paper is organized as follows. In Section 2, we
determine $\zeta_h(Q_n)=2^h(n-h)$ for $h\leq n-2$ and
$\eta_h(Q_n)=2^h(n-h)$ for $h\leq n-1$. In Section 3, we determine
$\zeta_{h}(S_{n})=\eta_{h}(S_{n})=h!(n-h)$ for $1\leq h \leq n-1$.
In Section 4, we determine $\eta_3(B_n)=6(n-3)$ for $n\geq 4$ and
point out a flaw in the proof of this conclusion in~\cite{yw14}. A
conclusion is in Section 5.

For graph terminology and notation not defined here we follow
Xu~\cite{x01}. For a subset $X$ of vertices in $G$, we do not
distinguish $X$ and the induced subgraph $G[X]$.

\section{Hypercubes}

\begin{figure}[h]
\psset{unit=0.9}
\begin{pspicture}(-3.3,-.3)(0,4)
\cnode(1,1){.1}{0}\rput(.75,1){\scriptsize0}
\cnode(1,3){.1}{1}\rput(.75,3){\scriptsize1}
\ncline{0}{1}\rput(1,.5){\scriptsize$Q_1$}
\end{pspicture}
\begin{pspicture}(-1.5,-.3)(2,3)
\cnode(1,1){.1}{00}\rput(.7,1){\scriptsize00}
\cnode(1,3){.1}{10}\rput(.7,3){\scriptsize10}
\cnode(3,1){.1}{01}\rput(3.35,1){\scriptsize01}
\cnode(3,3){.1}{11}\rput(3.3,3){\scriptsize11}
\ncline{00}{01}\ncline{01}{11}\ncline{11}{10}\ncline{10}{00}
\rput(2,.5){\scriptsize$Q_2$}
\end{pspicture}
\begin{pspicture}(-1.5,-.3)(2,3)
\cnode(1,1){.1}{000}\rput(.6,1){\scriptsize000}
\cnode(1,3){.1}{001}\rput(.6,3){\scriptsize001}
\cnode(3,1){.1}{100}\rput(3.4,1){\scriptsize100}
\cnode(3,3){.1}{101}\rput(3.4,3){\scriptsize101}
\cnode(1.7,1.7){.1}{010}\rput(2.1,2.){\scriptsize010}
\cnode(1.7,3.7){.1}{011}\rput(1.3,3.8){\scriptsize011}
\cnode(3.7,1.7){.1}{110}\rput(4.1,1.8){\scriptsize110}
\cnode(3.7,3.7){.1}{111}\rput(4.1,3.8){\scriptsize111}
\ncline{000}{001}\ncline{001}{101}\ncline{101}{100}\ncline{100}{000}
\ncline{010}{011}\ncline{011}{111}\ncline{111}{110}\ncline{110}{010}
\ncline{000}{010}\ncline{001}{011}\ncline{101}{111}\ncline{100}{110}
\rput(2,.5){\scriptsize$Q_3$}
\end{pspicture}
\vskip2pt
\begin{pspicture}(-4.2,0.5)(5,4)
\cnode(1,1){.1}{0000}\rput(.55,1){\scriptsize0000}
\cnode(1,3){.1}{0100}\rput(.55,3){\scriptsize0100}
\cnode(3,1){.1}{0001}\rput(2.62,1.2){\scriptsize0001}
\cnode(3,3){.1}{0101}\rput(2.62,3.2){\scriptsize0101}
\cnode(1.9,1.8){.1}{0010}\rput(1.45,1.85){\scriptsize0010}
\cnode(1.9,3.8){.1}{0110}\rput(1.45,3.85){\scriptsize0110}
\cnode(3.9,1.8){.1}{0011}\rput(3.51,2){\scriptsize0011}
\cnode(3.9,3.8){.1}{0111}\rput(3.51,4){\scriptsize0111}
\cnode(5,1){.1}{1001}\rput(5.4,.8){\scriptsize1001}
\cnode(5,3){.1}{1101}\rput(5.4,2.8){\scriptsize1101}
\cnode(7,1){.1}{1000}\rput(7.4,.85){\scriptsize1000}
\cnode(7,3){.1}{1100}\rput(7.4,2.85){\scriptsize1100}
\cnode(5.9,1.8){.1}{1011}\rput(6.3,2){\scriptsize1011}
\cnode(5.9,3.8){.1}{1111}\rput(6.3,4){\scriptsize1111}
\cnode(7.9,1.8){.1}{1010}\rput(8.34,1.8){\scriptsize1010}
\cnode(7.9,3.8){.1}{1110}\rput(8.34,3.8){\scriptsize1110}
\ncline{0000}{0001}\ncline{0001}{0101}\ncline{0101}{0100}\ncline{0100}{0000}
\ncline{0010}{0011}\ncline{0011}{0111}\ncline{0111}{0110}\ncline{0110}{0010}
\ncline{0000}{0010}\ncline{0001}{0011}\ncline{0101}{0111}\ncline{0100}{0110}
\ncline{1001}{1000}\ncline{1000}{1010}\ncline{1010}{1011}\ncline{1011}{1001}
\ncline{1101}{1100}\ncline{1100}{1110}\ncline{1110}{1111}\ncline{1111}{1101}
\ncline{1001}{1101}\ncline{1000}{1100}\ncline{1010}{1110}\ncline{1011}{1111}
\nccurve[angleA=-20,angleB=-160]{0000}{1000}
\nccurve[angleA=-20,angleB=-160]{0010}{1010}
\nccurve[angleA=-20,angleB=-160]{0001}{1001}
\nccurve[angleA=-20,angleB=-160]{0011}{1011}
\nccurve[angleA=20,angleB=160]{0100}{1100}
\nccurve[angleA=20,angleB=160]{0110}{1110}
\nccurve[angleA=20,angleB=160]{0101}{1101}
\nccurve[angleA=20,angleB=160]{0111}{1111}
\rput(1.5,.4){\scriptsize$Q_4$}
\end{pspicture}
\caption{
\label{f1}                                       
\footnotesize  The $n$-cubes $Q_1$, $Q_2$, $Q_3$ and $Q_4$}
\end{figure}
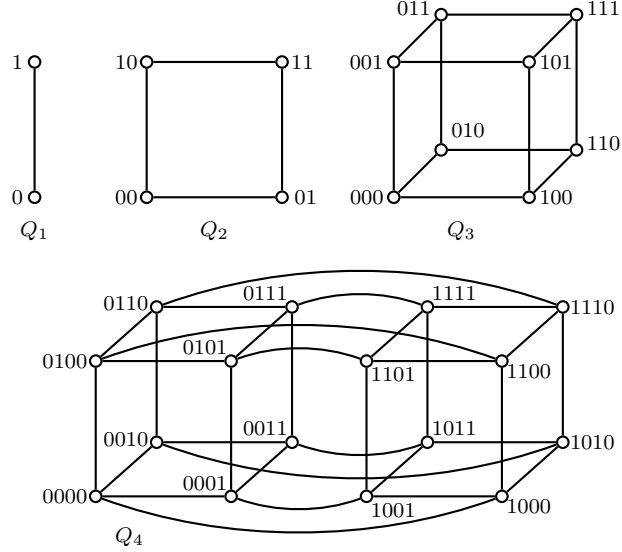

The hypercube $Q_n$ has the vertex-set consisting of $2^n$ binary
strings of length $n$,
two vertices being linked by an edge if and only if they differ in
exactly one position. Hypercubes $Q_2, Q_3, Q_4$ are shown in
Fig.~\ref{f1}.

The hypercube $Q_n$ is also defined as Cartesian product $K_2\times
K_2 \times \cdots \times K_2 $ of $n$ complete graph $K_2$. Thus,
$Q_n=Q_h\times Q_{n-h}$ for $1\leq h\leq n-1$, and $Q_n$ is a Cayley
graph with degree $n$ (see Xu~\cite{x01}). Oh {\it et
al.}~\cite{oc93} and Wu {\it et al.}~\cite{wg98} independently
determined $\kappa^{h}(Q_n)$, and Xu~\cite{x00c} determined
$\lambda^{h}(Q_n)$.

\begin{lem}\label{lem2.1}
$\kappa^h(Q_n)=2^h(n-h)$  for any $h$ with $0\leq h\leq n-2$, and
$\lambda^h(Q_n)=2^h(n-h)$ for any $h$ with $0\leq h\leq n-1$.
\end{lem}

Since $\delta_h=\delta(Q_h)=h$, by Lemma~\ref{lem1.1} and
Lemma~\ref{lem2.1}, the following corollary holds.

 \begin{cor}\label{cor2.2}
$\zeta_h(Q_n)\geq 2^h(n-h)$ for any $h$ with $1\leq h\leq n-2$, and
$\eta_h(Q_n)\geq 2^h(n-h)$ for any $h$ with $1\leq h\leq n-1$.
\end{cor}

 \begin{lem}\label{lem2.3}
$\zeta_h(Q_n)\leq 2^h(n-h)$ for any $h$ with $1\leq h\leq n-2$, and
$\eta_h(Q_n)\leq 2^h(n-h)$ for any $h$ with $1\leq h\leq n-1$.
\end{lem}

 \begin{pf}
For $h\leq n-2$, let $Q_n=Q_h\times Q_{n-h}$. Let $x$ be a vertex in
$Q_{n-h}$, $N(x)$ be the neighbor-set of $x$ in $Q_{n-h}$, and
$S=Q_h\times N(x)$. Then $Q_n-S$ is disconnected, and isomorphic to
$Q_h\times (Q_{n-h}-N(x))$, so each vertex of $Q_n-S$ is in some
$Q_h$. It follows that $\zeta_h(Q_n) \leq |S| = 2^h(n-h)$.

For $h\leq n-1$, $Q_n=Q_h\times Q_{n-h}$, let $F$ be the set of
edges between $Q_h$  and $Q_n-Q_h$. It is easy to see that $Q_n-F$
is disconnected, and isomorphic to $Q_h\times (Q_{n-h}-E(y))$, where
$E(y)$ is the set of edges in $Q_{n-h}$ incident with a vertex $y$
of $Q_{n-h}$, Thus each vertex of $Q_n-F$ is in some $Q_h$. It
follows that $\eta_h(Q_n) \leq |F| = 2^h(n-h)$.
 \end{pf}

\vskip6pt

Combining Corollary~\ref{cor2.2} with Lemma~\ref{lem2.3}, we obtain
the following conclusion.

 \begin{thm}\label{thm2.4}
$\zeta_h(Q_n)=2^h(n-h)$ for any $h$ with $1\leq h\leq n-2$, and
$\eta_h(Q_n)=2^h(n-h)$ for any $h$ with $1\leq h\leq n-1$.
\end{thm}

\section{Star graphs}

For a given integer $n$ with $n\geq 2$, let $I_n=\{1,2,\ldots,n\}$,
$I'_n=\{2,\ldots,n\}$ and $P(n)=\{ p_{1}p_{2}\ldots p_{n}:\ p_{i}\in
I_n, p_{i}\neq p_{j}, 1\leq i\neq j\leq n\}$, the set of
permutations on $I_n$. Clearly, $|P(n)|=n\,!$. For $p=p_{1}\ldots
p_j\ldots p_{n}\in P(n)$, the digit $p_j$ is called the $j$-th digit
of $p$.

\begin{figure}[h]
\begin{center}
\psset{unit=0.9}
\begin{pspicture}(1,-2.5)(4,4.)
\rput{0}{%
\SpecialCoor\degrees[6]
 \ncline{10}{11}\ncline{11}{12}
 \multido{\i=0+1}{6}{\cnode(.9;\i){.11}{1\i}}
 \ncline{10}{11}\ncline{11}{12}\ncline{12}{13}\ncline{13}{14}\ncline{14}{15}\ncline{15}{10}
 }%
 \cnode(-.5,3.5){.11}{a1}\rput(-.9,3.5){\scriptsize$12$}
 \cnode(.5,3.5){.11}{a2}\rput(.9,3.5){\scriptsize$21$}
 \ncline{a1}{a2}\rput(0,3.){\scriptsize$S_2$}\rput(0,-1.2){\scriptsize$S_3$}
 \rput(-.8,0.8){\scriptsize$123$}\rput(-1.25,-0.){\scriptsize$321$}\rput(-.8,-.8){\scriptsize$231$}
 \rput(.85,0.8){\scriptsize$213$}\rput(1.3,-0.){\scriptsize$312$}\rput(.8,-.8){\scriptsize$132$}
 \end{pspicture}
\begin{pspicture}(-2,-3.7)(2,4.2)
\rput{0}{%
\SpecialCoor\degrees[6]
 \multido{\i=0+1}{6}{\cnode(.9;\i){.11}{1\i}}
 \ncline{10}{11}\ncline{11}{12}\ncline{12}{13}\ncline{13}{14}\ncline{14}{15}\ncline{15}{10}
 \multido{\i=0+1}{6}{\cnode(1.8;\i){.11}{2\i}}
 \ncline{20}{21}\ncline{21}{22}\ncline{22}{23}\ncline{23}{24}\ncline{24}{25}\ncline{25}{20}
 \multido{\i=0+1}{6}{\cnode(2.8;\i){.11}{3\i}}
 \ncline{30}{31}\ncline{31}{32}\ncline{32}{33}\ncline{33}{34}\ncline{34}{35}\ncline{35}{30}
 }%
 \rput{90}{%
\SpecialCoor\degrees[6]
 \multido{\i=0+1}{6}{\cnode(3.5;\i){.11}{4\i}}
}%
 \ncline[linecolor=blue,linewidth=1pt]{40}{13}\ncline[linecolor=blue,linewidth=1pt]{40}{22}\ncline[linecolor=blue,linewidth=1pt]{40}{31}
 \ncline[linecolor=blue,linewidth=1pt]{41}{12}\ncline[linecolor=blue,linewidth=1pt]{41}{23}\ncline[linecolor=blue,linewidth=1pt]{41}{34}
 \ncline[linecolor=blue,linewidth=1pt]{42}{15}\ncline[linecolor=blue,linewidth=1pt]{42}{24}\ncline[linecolor=blue,linewidth=1pt]{42}{33}
 \ncline[linecolor=blue,linewidth=1pt]{43}{14}\ncline[linecolor=blue,linewidth=1pt]{43}{25}\ncline[linecolor=blue,linewidth=1pt]{43}{30}
 \ncline[linecolor=blue,linewidth=1pt]{44}{11}\ncline[linecolor=blue,linewidth=1pt]{44}{20}\ncline[linecolor=blue,linewidth=1pt]{44}{35}
 \ncline[linecolor=blue,linewidth=1pt]{45}{10}\ncline[linecolor=blue,linewidth=1pt]{45}{21}\ncline[linecolor=blue,linewidth=1pt]{45}{32}

 \rput(0,3.8){\scriptsize$1324$}\rput(0,-3.8){\scriptsize$1423$}
 \rput(-1.5,2.7){\scriptsize$2341$}\rput(1.7,2.7){\scriptsize$4321$}
 \rput(-1.5,-2.7){\scriptsize$4231$}\rput(1.7,-2.7){\scriptsize$2431$}
 \rput(-3.5,1.7){\scriptsize$1234$}\rput(3.5,1.7){\scriptsize$1342$}
 \rput(-3.5,-1.7){\scriptsize$1243$}\rput(3.5,-1.7){\scriptsize$1432$}
 \rput(-3.3,0.){\scriptsize$3241$}\rput(3.3,0.){\scriptsize$3421$}
 \rput(-1.35,0.){\scriptsize$3214$}\rput(2.25,0.){\scriptsize$3412$}
 \rput(-.4,0.){\scriptsize$3124$}\rput(.45,0.){\scriptsize$3142$}
 \rput(-1.3,1.7){\scriptsize$2314$}\rput(.8,1.85){\scriptsize$4312$}
 \rput(-1.,-1.8){\scriptsize$4213$}\rput(.6,-1.3){\scriptsize$2413$}
 \rput(-.2,0.5){\scriptsize$2134$}\rput(.4,1.05){\scriptsize$4132$}
 \rput(-.2,-0.5){\scriptsize$4123$}\rput(.85,-.85){\scriptsize$2143$}

 \rput(-1.5,-3.7){\scriptsize$S_4$}

 \end{pspicture}

\caption{\label{f2}\footnotesize {The star graphs $S_2, S_3$ and
$S_4$} }
\end{center}
\end{figure}
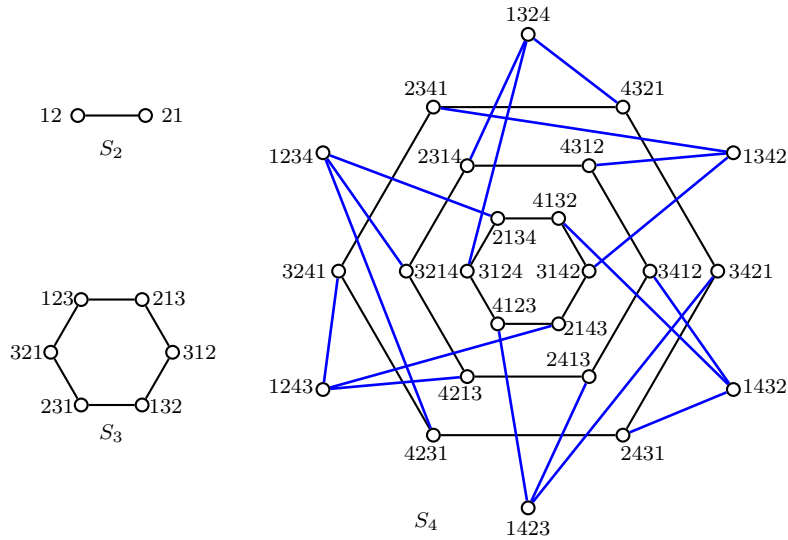

The $n$-dimensional star graph $S_{n}$ has vertex-set $P(n)$ and has
an edge between any two vertices if and only if one can be obtained
from the other by swapping the $1$-th digit and the $i$-th digit for
$i\in I'_n$, that is, two vertices $x=p_{1}p_{2}\ldots p_{i}\ldots
p_{n}$ and $y$ are adjacent if and only if $y=p_{i}p_{2}\ldots
p_{i-1}p_1p_{i+1}\ldots p_{n}$ for some $i\in I'_n$. The star graphs
$S_2, S_3, S_4$ are shown in Fig.~\ref{f2}. It is shown that the
star graph $S_n$ is a Cayley graph with degree $(n-1)$ (see Akers
and Krishnamurthy~\cite{ak89}).

For a fixed $i\in I_n$, let $S^{j:i}_n$ denote a subgraph of $S_{n}$
induced by all vertices whose the $j$-th digit is $i$ for each $j\in
I_n$. By definition, it is easy to see that $S^{j:i}_n$ is
isomorphic to $S_{n-1}$ for each $j\in I'_n$ and $S^{1:i}_n$ is an
edgeless graph with $(n-1)\,!$ vertices. As shown in Fig.~\ref{f1},
$S^{1:1}_4$ is an edgeless graph with $6$ vertices, $S^{j:1}_4$ is
isomorphic to $S_3$ for each $j$ with $2\leq j\leq 4$.

There are two different hierarchical structures of $S_n$ depending
on different partition methods. The first one is partitioned along a
fixed dimension, which is well-known and used frequently. The second
one is partitioned along a fixed digit in $I_n$, which is a new
structure proposed by Shi~{\it et al.}\cite{sls12}.

\begin{lem}\label{lem3.2}\textnormal{(Shi {\it et al.}~\cite{sls12}, 2012)}
For a fixed $i\in I_n$, $S_{n}$ can be partitioned into $n$
subgraphs $S^{j:i}_n$, which is isomorphic to $S_{n-1}$ for each
$j\in I'_n$ and $S^{1:i}_n$ is an edgeless graph with $(n - 1)\,!$
vertices. Moreover, there are $(n-2)!$ independent edges between
$S^{1:i}_n$ and $S^{j:i}_n$ for any $j\in I'_n$, and there is no
edge between $S^{j_1:i}_n$ and $S^{j_2:i}_n$ for any $j_1,j_2\in
I'_n$ with $j_1\ne j_2$.
\end{lem}

Clearly, $S_1,S_2,S_3$ are isomorphic to $K_1,K_2,C_6$,
respectively. As shown in Fig.~\ref{f1}, $S_4$ is partitioned along
digit 1.

\begin{lem}\label{lem3.3}\textnormal{(Li and Xu~\cite{lx14}, 2014)}
$\kappa^{h}(S_{n})=\lambda^h(S_{n}) =(h+1)\,!\,(n-h-1)$ for any $h$
with $0\leq h\leq n-2$.
\end{lem}

Since $\delta_h=\delta(S_h)=h-1$, by Lemma~\ref{lem1.1} and
Lemma~\ref{lem3.3}, the following corollary holds.

\begin{cor}\label{cor3.3}
$\zeta_h(S_{n})\geq h\,!\,(n-h)$ and $\eta_h(S_{n})\geq h\,!\,(n-h)$
for any $h$ with $1\leq h \leq n-1$.
\end{cor}

To determine $\zeta_h(S_{n})$ and $\eta_h(S_{n})$, we investigate
their upper bounds.

\begin{lem}\label{lem3.4}\textnormal{(Yang {\it et al.}~\cite{yw12})}
$\eta_h(S_{n})\leq h\,!\,(n-h)$ for any $h$ with $1\leq h \leq n-1$.
\end{lem}

Now we establish the upper bound on $\zeta_h(S_{n})$ by
Lemma~\ref{lem3.2}.

\begin{lem}\label{lem3.5}
$\zeta_h(S_{n})\leq h\,!\,(n-h)$ for any $h$ with $1\leq h \leq
n-1$.
\end{lem}

\begin{pf}
Let
 $$
 X=\{ \ p_1p_2 \cdots  p_{h}12\cdots (n-h)\in V(S_n):
 \  p_i \in I_n \setminus I_{n-h} , \ i \ \in  I_{h}\}.
 $$
Then, $S_n[X]\cong S_{h}$. Let $T$ be the neighbor-set of $X$ in
$S_n-X$. By the definition of $S_n$, for a vertex of $X$, since it
has $h-1$ neighbors in $X$, it exactly has $(n-h)$ neighbors in $T$,
and every vertex of $T$ exactly has one neighbor in $X$. It follows
that
 $$
  |T|=h\,!\,(n-h).
 $$

Next we show $T$ is an $h$-embedded restricted vertex-cut of $S_n$.
It suffice to show each vertex on $S_n-(X\cup T)$ is in some
subgraph $S_h$ of $S_n-(X\cup T)$.

Assume that $u=p'_1p'_2 \cdots  p'_{n}$ is a vertex in $S_n-(X\cup T
)$, and let
 $$
 J=\{j\in I_{n-h}:\ p'_{h+j}\ne j \}\ {\rm and}\ J'=\{j\in J:\ p'_1\ne j\}.
 $$
Since $u\not \in X$, we have $J\not \ne \emptyset$. We claim $J'\ne
\emptyset$. Suppose to the contrary $J'=\emptyset$. Then $p'_1=j$
for each $j\in J$, and so $|J|$=1. Assume $J=\{j\}$. Note that
$p'_1=j$ and $p'_{h+j}\ne j$ $(1\leq j\leq n-h)$. Thus $u$ is a
neighbor of some vertex in $X$, that is, $u\in T$, which contradicts
to $u \notin T$.

Thus, $J'\ne \emptyset$. Let $j_0\in J'$ and $p'_{i_0}=j_0$ $(1\leq
i_0\leq n)$. Then $h+j_0\ne i_0$ and $p'_1\ne j_0$. We partition
$S_n$ by fixing digit $j_0$, Then $X\subseteq S^{(h+j_0):j_0}_{n}$
and $u\in S^{i_0:j_0}_{n}$. By Lemma~\ref{lem3.2},
$S^{(h+j_0):j_0}_{n}$ and $S^{i_0:j_0}_{n}$ are both isomorphic to
$S_{n-1}$, and there is no edge between $S^{(h+j_0):j_0}_{n}$ and
$S^{i_0:j_0}_{n}$, and so $u$ is in some $S_{n-1}$ of $S_n\setminus
(X\cup T )$, which implies that $u$ is in some $S_{h}$ of
$S_n\setminus (X\cup T )$. Thus, $T$ is an $h$-embedded restricted
vertex-cut of $S_n$, and so
 $$
 \zeta^{(h)}(S_n)\leq |T|=h\,!(n-h).
 $$
The lemma follows.\end{pf}

\begin{thm}\label{thm3.6}
$\zeta_{h}(S_{n})=\eta_{h}(S_{n})=h\,!\,(n-h)$ for any $h$ with
$1\leq h \leq n-1$.
\end{thm}

\begin{pf}
For $1\leq h \leq n-1$, combining Corollary~\ref{cor3.3} with
Lemma~\ref{lem3.4}, we have $\eta_{h}(S_{n})=h!(n-h)$, and combining
Corollary~\ref{cor3.3} with Lemma~\ref{lem3.5}, we have
$\zeta_{h}(S_{n})=h\,!(n-h)$.
\end{pf}


\section{Bubble-sort graphs}

The  $n$-dimensional bubble-sort graph $B_n$ has vertex-set $P(n)$
and has an edge between any two vertices if and only if one can be
obtained from the other by swapping the $i$-th digit and the
$(i+1)$-th digit where $1\leq i\leq n-1$. The bubble-sort graphs
$B_2,B_3$ and $B_4$ are shown in Fig.~\ref{f2}.

\begin{figure}[h]
\begin{center}
\psset{unit=0.9pt}
\begin{pspicture}(0,10)(280,170)

\cnode(0,140){3}{a} \cnode(30,140){3}{b} \ncline{a}{b}

\rput(-12,140){{\tiny $12$}} \rput(42,140){{\tiny $21$}}

\rput(15,130){{\scriptsize $B_2$}}


\cnode(15,30){3}{a1} \cnode(15,78){3}{a4} \cnode(-6,42){3}{a6}
\cnode(-6,66){3}{a5} \cnode(36,42){3}{a2} \cnode(36,66){3}{a3}
\ncline{a1}{a2} \ncline{a2}{a3} \ncline{a3}{a4} \ncline{a4}{a5}
\ncline{a5}{a6} \ncline{a6}{a1}

\rput(48,40){{\tiny $231$}} \rput(48,68){{\tiny $213$}}
\rput(-17,40){{\tiny $312$}} \rput(-17,68){{\tiny $132$}}
\rput(15,87){{\tiny $123$}} \rput(15,21){{\tiny $321$}}

\rput(15,5){{\scriptsize $B_3$}}



\cnode(120,74){3}{b6} \cnode(120,98){3}{b5} \cnode(141,110){3}{b4}
\cnode(141,62){3}{b1} \cnode(162,98){3}{b3} \cnode(162,74){3}{b2}
\ncline{b1}{b2} \ncline{b2}{b3} \ncline{b3}{b4} \ncline{b4}{b5}
\ncline{b5}{b6} \ncline{b6}{b1}

\cnode(186,74){3}{c6} \cnode(186,98){3}{c5} \cnode(207,110){3}{c4}
\cnode(207,62){3}{c1} \cnode(228,98){3}{c3} \cnode(228,74){3}{c2}
\ncline{c1}{c2} \ncline{c2}{c3} \ncline{c3}{c4} \ncline{c4}{c5}
\ncline{c5}{c6} \ncline{c6}{c1}

\cnode(219,17){3}{d6} \cnode(219,41){3}{d5} \cnode(240,53){3}{d4}
\cnode(240,5){3}{d1} \cnode(261,41){3}{d3} \cnode(261,17){3}{d2}
\ncline{d1}{d2} \ncline{d2}{d3} \ncline{d3}{d4} \ncline{d4}{d5}
\ncline{d5}{d6} \ncline{d6}{d1}

\cnode(219,131){3}{e6} \cnode(219,155){3}{e5} \cnode(240,167){3}{e4}
\cnode(240,119){3}{e1} \cnode(261,155){3}{e3} \cnode(261,131){3}{e2}
\ncline{e1}{e2} \ncline{e2}{e3} \ncline{e3}{e4} \ncline{e4}{e5}
\ncline{e5}{e6} \ncline{e6}{e1}

\ncline{b2}{c6} \ncline{b3}{c5} \ncline{d5}{c1} \ncline{d4}{c2}
\ncline{e1}{c3} \ncline{e6}{c4}

\ncline{b1}{d6} \ncline{e2}{d3} \ncline{e5}{b4}


\psarc(207,86){88}{-50}{50} \psarc(207,86){88}{70}{170}
\psarc(207,86){88}{190}{290}


\rput(140,5){{\scriptsize $B_4$}}

\rput(107,74){{\tiny $4312$}} \rput(107,98){{\tiny $4132$}}
\rput(141,119){{\tiny $1432$}} \rput(141,53){{\tiny $3412$}}
\rput(162,108){{\tiny $1342$}} \rput(162,64){{\tiny $3142$}}

\rput(186,64){{\tiny $3124$}} \rput(186,108){{\tiny $1324$}}
\rput(200,120){{\tiny $1234$}} \rput(200,52){{\tiny $3214$}}
\rput(242,98){{\tiny $2134$}} \rput(242,74){{\tiny $2314$}}

\rput(214,9){{\tiny $3421$}} \rput(208,40){{\tiny $3241$}}
\rput(252,56){{\tiny $2341$}} \rput(250,0){{\tiny $4321$}}
\rput(273,44){{\tiny $2431$}} \rput(275,17){{\tiny $4231$}}

\rput(207,133){{\tiny $1243$}} \rput(215,163){{\tiny $1423$}}
\rput(240,177){{\tiny $4123$}} \rput(252,117){{\tiny $2143$}}
\rput(273,159){{\tiny $4213$}} \rput(273,128){{\tiny $2413$}}

\end{pspicture}

\end{center}
\caption{\label{f3}\footnotesize{The bubble-sort graphs $B_2$, $B_3$
and $B_4$. }}
\end{figure}
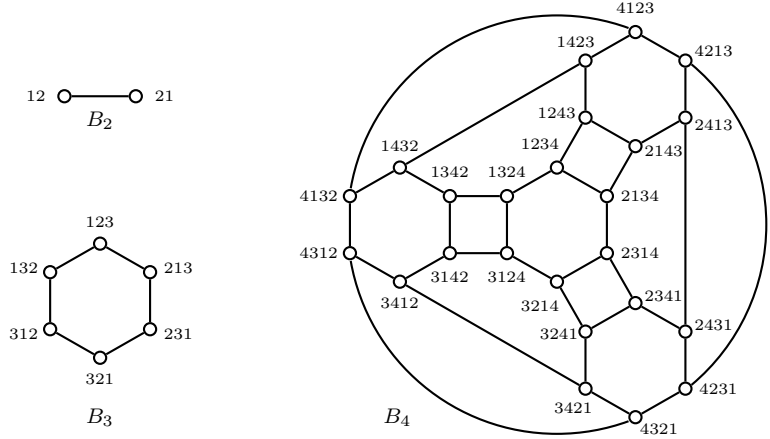

It is shown that the bubble-sort graph $B_n$ is a Cayley graph with
degree $(n-1)$ (see Akers and Krishnamurthy~\cite{ak89}). More
specifically, $B_n$ is bipartite and contains $n$ disjoint
sub-graphs $B_{n-1}$ by fixing the $1$-th digit or $n$-th digit (see
Fig.~\ref{f3}).

For each $i, t\in I_n$, let $B^{t:i}_{n}$ denote a subgraph of
$B_{n}$ induced by all vertices whose the $t$-th digit is $i$.
Clearly, $B^{t:i}_{n}\cong B_{n-1}$ for each $i, t\in I_n$.

\begin{lem}\label{lem4.1}\textnormal{(Akers and Krishnamurthy~\cite{ak89}, 1989)}
For a fixed $t \in \{1, n\}$, $B_{n}$ can be partitioned into $n$
subgraphs $B^{t:j}_n$ isomorphic to $B_{n-1}$ for each $j\in I_n$,
moreover, there are $(n-2)!$ independent edges between $B^{t:j_1}_n$
and $B^{t:j_2}_n$ for any $j_1,j_2\in I_n$ with $j_1\ne j_2$.
\end{lem}

\begin{lem}\label{lem4.2}\textnormal{(Yang {\it et al.}~\cite{yw14}, 2014)}\
$\eta_{h}(B_{n})\leq h\,!\,(n-h)$ for any $h$ with $1\leq h\leq n-1$
and $n\geq 2$.
\end{lem}

It was showed by Xu~\cite{x00d} in 2000 that for a vertex-transitive
connected graph $G$ with order $n\ (\geqslant 4)$ and with degree
$d\ (\geqslant 2)$, $\lambda^1(G)=2d-2$ if $n$ is odd or $G$
contains no triangles. Since $B_n$ is vertex-transitive and
bipartite, we have the following result.

\begin{lem}\label{lem4.3}
$\lambda^{1}(B_{n})=2n-4$ for $n\geq 3$.
\end{lem}

Since $\delta_h=\delta(B_h)=h-1$, by Lemma~\ref{lem1.1} and
Lemma~\ref{lem4.3} $\eta_2(B_n)\geq\lambda^1(B_n)=2n-4$, and by
Lemma~\ref{lem4.2} the following corollary holds.

\begin{cor}\
$\eta_{2}(B_{n})=2n-4$ if $n\geq 3$.
\end{cor}

Yang {\it et al.} also stated this result (see Lemma 4.4
in~\cite{yw14}), however, their proof is not true since they use an
incorrect condition that the girth of $B_n$ is 6 (in fact the girth
of $B_n$ is 4 for $n\geq4$).

For $h=3$, using Lemma~\ref{lem4.2} and Lemma~\ref{lem1.1}, Yang
{\it et al.}~\cite{yw14} gave $\eta_{3}(B_{n})=6(n-3)$. However,
when Lemma~\ref{lem1.1} is used, they assume an incorrect
conclusion: $\lambda^2(B_n)=6(n-3)$ (see Lemma 4.6 in~\cite{yw14}).
In fact, $\lambda^2(B_n)=4(n-3)$ for $n\geq 4$ (see \cite{ylm10}).
It follows that Yang {\it et al.}'s proof is not correct.

In this section, we will prove $\eta_3(B_n)\geq 6(n-3)$ for $n\geq
4$. The proof proceeds by induction on $n\, (\geq 4)$. And so we
first consider the case of $n=4$.

An edge-cut $F$ of $B_n$ is called a $B_h$-edge-cut if every
component of $B_n-F$ contains a $B_h$ as subgraph. By definition,
every $h$-embedded vertex-cut of $B_n$ is certainly a
$B_h$-edge-cut.

\begin{lem}\label{lem4.5}
If $F$ be a $B_3$-edge-cut of $B_4$, then $|F|\geq 6$.
\end{lem}

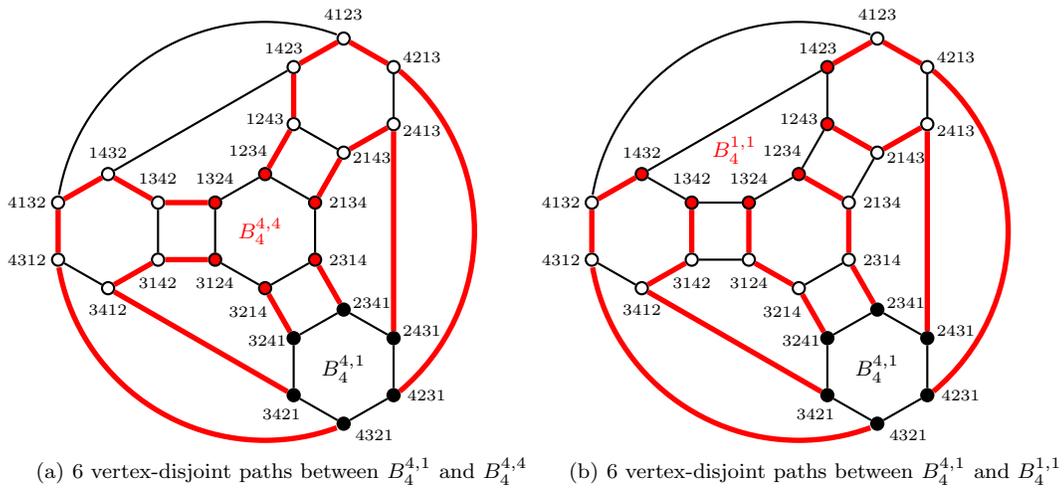
\begin{figure}[h]
\begin{center}
\psset{unit=0.9pt}
\begin{pspicture}(280,0)(310,190)
\cnode(120,74){3}{b6} \cnode(120,98){3}{b5} \cnode(141,110){3}{b4}
\cnode(141,62){3}{b1} \cnode(162,98){3}{b3} \cnode(162,74){3}{b2}
\ncline[linecolor=red,linewidth=2pt]{b1}{b2} \ncline{b2}{b3}
\ncline[linecolor=red,linewidth=2pt]{b3}{b4}
\ncline[linecolor=red,linewidth=2pt]{b4}{b5}
\ncline[linecolor=red,linewidth=2pt]{b5}{b6} \ncline{b6}{b1}

\cnode[fillstyle=solid,fillcolor=red](186,74){3}{c6}
\cnode[fillstyle=solid,fillcolor=red](186,98){3}{c5}
\cnode[fillstyle=solid,fillcolor=red](207,110){3}{c4}
\cnode[fillstyle=solid,fillcolor=red](207,62){3}{c1}
\cnode[fillstyle=solid,fillcolor=red](228,98){3}{c3}
\cnode[fillstyle=solid,fillcolor=red](228,74){3}{c2}
 \ncline{c1}{c2}
\ncline{c2}{c3} \ncline{c3}{c4} \ncline{c4}{c5} \ncline{c5}{c6}
\ncline{c6}{c1}

\cnode*(219,17){3}{d6} \cnode*(219,41){3}{d5} \cnode*(240,53){3}{d4}
\cnode*(240,5){3}{d1} \cnode*(261,41){3}{d3} \cnode*(261,17){3}{d2}
\ncline{d1}{d2} \ncline{d2}{d3} \ncline{d3}{d4} \ncline{d4}{d5}
\ncline{d5}{d6} \ncline{d6}{d1}

\cnode(219,131){3}{e6} \cnode(219,155){3}{e5} \cnode(240,167){3}{e4}
\cnode(240,119){3}{e1} \cnode(261,155){3}{e3} \cnode(261,131){3}{e2}
\ncline[linecolor=red,linewidth=2pt]{e1}{e2} \ncline{e2}{e3}
\ncline[linecolor=red,linewidth=2pt]{e3}{e4}
\ncline[linecolor=red,linewidth=2pt]{e4}{e5}
\ncline[linecolor=red,linewidth=2pt]{e5}{e6} \ncline{e6}{e1}

\ncline[linecolor=red,linewidth=2pt]{b2}{c6}
\ncline[linecolor=red,linewidth=2pt]{b3}{c5}
\ncline[linecolor=red,linewidth=2pt]{d5}{c1}
\ncline[linecolor=red,linewidth=2pt]{d4}{c2}
\ncline[linecolor=red,linewidth=2pt]{e1}{c3}
\ncline[linecolor=red,linewidth=2pt]{e6}{c4}

\ncline[linecolor=red,linewidth=2pt]{b1}{d6}
\ncline[linecolor=red,linewidth=2pt]{e2}{d3} \ncline{e5}{b4}


\psarc[linecolor=red,linewidth=2pt](207,86){88}{-50}{50}
\psarc(207,86){88}{70}{170}
\psarc[linecolor=red,linewidth=2pt](207,86){88}{190}{290}


 \rput(240,28){\scriptsize $B_4^{4,1}$}
 \rput(205,86){\scriptsize {\red $B_4^{4,4}$}}

\rput(107,74){{\tiny $4312$}} \rput(107,98){{\tiny $4132$}}
\rput(141,119){{\tiny $1432$}} \rput(141,53){{\tiny $3412$}}
\rput(162,108){{\tiny $1342$}} \rput(162,64){{\tiny $3142$}}

\rput(186,64){{\tiny $3124$}} \rput(186,108){{\tiny $1324$}}
\rput(200,120){{\tiny $1234$}} \rput(200,52){{\tiny $3214$}}
\rput(242,98){{\tiny $2134$}} \rput(242,74){{\tiny $2314$}}

\rput(214,9){{\tiny $3421$}} \rput(208,40){{\tiny $3241$}}
\rput(252,56){{\tiny $2341$}} \rput(253,0){{\tiny $4321$}}
\rput(273,44){{\tiny $2431$}} \rput(275,17){{\tiny $4231$}}

\rput(207,133){{\tiny $1243$}} \rput(215,163){{\tiny $1423$}}
\rput(240,177){{\tiny $4123$}} \rput(252,117){{\tiny $2143$}}
\rput(273,159){{\tiny $4213$}} \rput(273,128){{\tiny $2413$}}
 \rput(214,-15){\scriptsize (a) 6 vertex-disjoint paths between $B_4^{4,1}$ and $B_4^{4,4}$}

\end{pspicture}
\begin{pspicture}(90,0)(120,180)
\cnode(120,74){3}{b6} \cnode(120,98){3}{b5}
\cnode[fillstyle=solid,fillcolor=red](141,110){3}{b4}
\cnode(141,62){3}{b1}
\cnode[fillstyle=solid,fillcolor=red](162,98){3}{b3}
\cnode(162,74){3}{b2} \ncline[linecolor=red,linewidth=2pt]{b1}{b2}
\ncline[linecolor=red,linewidth=2pt]{b2}{b3} \ncline{b3}{b4}
\ncline[linecolor=red,linewidth=2pt]{b4}{b5}
\ncline[linecolor=red,linewidth=2pt]{b5}{b6} \ncline{b6}{b1}

\cnode(186,74){3}{c6}
\cnode[fillstyle=solid,fillcolor=red](186,98){3}{c5}
\cnode[fillstyle=solid,fillcolor=red](207,110){3}{c4}
\cnode(207,62){3}{c1} \cnode(228,98){3}{c3} \cnode(228,74){3}{c2}
\ncline{c1}{c2} \ncline[linecolor=red,linewidth=2pt]{c2}{c3}
\ncline[linecolor=red,linewidth=2pt]{c3}{c4} \ncline{c4}{c5}
\ncline[linecolor=red,linewidth=2pt]{c5}{c6}
\ncline[linecolor=red,linewidth=2pt]{c6}{c1}

\cnode*(219,17){3}{d6} \cnode*(219,41){3}{d5} \cnode*(240,53){3}{d4}
\cnode*(240,5){3}{d1} \cnode*(261,41){3}{d3} \cnode*(261,17){3}{d2}
\ncline{d1}{d2} \ncline{d2}{d3} \ncline{d3}{d4} \ncline{d4}{d5}
\ncline{d5}{d6} \ncline{d6}{d1}

\cnode[fillstyle=solid,fillcolor=red](219,131){3}{e6}
\cnode[fillstyle=solid,fillcolor=red](219,155){3}{e5}
\cnode(240,167){3}{e4} \cnode(240,119){3}{e1} \cnode(261,155){3}{e3}
\cnode(261,131){3}{e2} \ncline[linecolor=red,linewidth=2pt]{e1}{e2}
\ncline{e2}{e3} \ncline[linecolor=red,linewidth=2pt]{e3}{e4}
\ncline[linecolor=red,linewidth=2pt]{e4}{e5} \ncline{e5}{e6}
\ncline[linecolor=red,linewidth=2pt]{e6}{e1}

\ncline{b2}{c6} \ncline{b3}{c5}
\ncline[linecolor=red,linewidth=2pt]{d5}{c1}
\ncline[linecolor=red,linewidth=2pt]{d4}{c2} \ncline{e1}{c3}
\ncline{e6}{c4}

\ncline[linecolor=red,linewidth=2pt]{b1}{d6}
\ncline[linecolor=red,linewidth=2pt]{e2}{d3} \ncline{e5}{b4}


\psarc[linecolor=red,linewidth=2pt](207,86){88}{-50}{50}
\psarc(207,86){88}{70}{170}
\psarc[linecolor=red,linewidth=2pt](207,86){88}{190}{290}


 \rput(240,28){\scriptsize $B_4^{4,1}$}
 \rput(180,120){\scriptsize {\red $B_4^{1,1}$}}

\rput(107,74){{\tiny $4312$}} \rput(107,98){{\tiny $4132$}}
\rput(141,119){{\tiny $1432$}} \rput(141,53){{\tiny $3412$}}
\rput(162,108){{\tiny $1342$}} \rput(162,64){{\tiny $3142$}}

\rput(186,64){{\tiny $3124$}} \rput(186,108){{\tiny $1324$}}
\rput(200,120){{\tiny $1234$}} \rput(200,52){{\tiny $3214$}}
\rput(242,98){{\tiny $2134$}} \rput(242,74){{\tiny $2314$}}

\rput(214,9){{\tiny $3421$}} \rput(208,40){{\tiny $3241$}}
\rput(252,56){{\tiny $2341$}} \rput(253,0){{\tiny $4321$}}
\rput(273,44){{\tiny $2431$}} \rput(275,17){{\tiny $4231$}}

\rput(207,133){{\tiny $1243$}} \rput(215,163){{\tiny $1423$}}
\rput(240,177){{\tiny $4123$}} \rput(252,117){{\tiny $2143$}}
\rput(273,159){{\tiny $4213$}} \rput(273,128){{\tiny $2413$}}

\rput(214,-15){\scriptsize (b) 6 vertex-disjoint paths between
$B_4^{4,1}$ and $B_4^{1,1}$}

\end{pspicture}

\end{center}
\caption{\label{f4}\footnotesize{The disjoint paths between two
disjoint $B_3$s in $B_4$ }}
\end{figure}

\begin{pf}
We first prove that there are 6 vertex-disjoint paths between each
two disjoint $B_3$s in $B_4$. By symmetry, we only need to consider
the paths from $B_4^{4,1}$ to $B_4^{4,4}$ and $B_4^{1,1}$. Such
paths are illustrated in the Table~\ref{tb1} (also see
Fig.~\ref{f4}). If $F$ is a $B_3$-edge-cut of $B_4$, then $F$
contains at least one edge of every path in the 6 vertex-disjoint
paths between the two disjoint $B_3$s. It follows that $|F|\geq 6$.
 \end{pf}

\begin{table}[h]
\begin{tabular}
{|p{53pt}|p{200pt}|p{150pt}|}
\hline
From  $\backslash$ {To}&
  {\footnotesize \hspace{2.5cm}  $B_4^{4,4}$}&
  {\footnotesize \hspace{1.7cm} $B_4^{1,1}$} \\
\hline \multirow{1}{*}[-3.2em]{\footnotesize
\hspace{0.5cm}$B_4^{4,1}$}& {\small 324}\textbf{\small 1} - {\small
321}\textbf{\small 4}
\par {\small 234}\textbf{\small 1} - {\small 231}\textbf{\small 4}
\par {\small 342}\textbf{\small 1} - {\small 3412} - {\small 3142} - {\small 312}\textbf{\small 4}
\par {\small 243}\textbf{\small 1} - {\small 2413} - {\small 2143} - {\small 213}\textbf{\small 4}
 \par {\small 432}\textbf{\small 1} - {\small 4312} - {\small 4132} - {\small 1432} - {\small 1342} - {\small 132}\textbf{\small 4}
 \par {\small 423}\textbf{\small 1} - {\small 4213} - {\small 4123} - {\small 1423} -
 {\small 1243} - {\small 123}\textbf{\small 4} &
 {\small 324}\textbf{\small 1} - {\small 3214} - {\small 3124} - \textbf{\small 1}{\small 324}
\par {\small 234}\textbf{\small 1} - {\small 2314} - {\small 2134} - \textbf{\small 1}{\small 234}
 \par {\small 342}\textbf{\small 1} - {\small 3412} - {\small 3142} - \textbf{\small 1}{\small 342}
\par {\small 243}\textbf{\small 1} - {\small 2413} - {\small 2143} - \textbf{\small 1}{\small 243}
\par {\small 432}\textbf{\small 1} - {\small 4312} - {\small 4132} - \textbf{\small 1}{\small 432}
 \par {\small 423}\textbf{\small 1} - {\small 4213} - {\small 4123} - \textbf{\small 1}{\small 423} \\
\hline
\end{tabular}
 \caption{The disjoint paths between two disjoint $B_3$s in $B_4$\label{tb1}}
\end{table}

\begin{thm}\label{thm4.6}
$\eta_{3}(B_{n})= 6(n-3)$ if $n\geq 4$.
\end{thm}

\begin{pf}
By Lemma~\ref{lem4.2}, we only need to prove $\eta_{3}(B_{n})\geq
6(n-3)$ for $n\geq 4$. Let $F$ be a $B_3$-edge-cut of $B_n$. Since
any $3$-embedded edge-cut is certainly a $B_3$-edge-cut of $B_n$, we
have that
  $
  \eta_{3}(B_{n})\geq |F|.
  $
Thus it suffice to show
\begin{equation}\label{e4.1}
  |F|\geq 6(n-3).
  \end{equation}
The proof proceeds by induction on $n\, (\geq 4)$. When $n=4$, the
result follows by Lemma~\ref{lem4.5}. Assume the induction
hypothesis for any $m$ with $4\leq m\leq n-1$ with $n\geq 5$, that
is
\begin{equation}\label{e4.2}
  |F'|\geq 6(m-3)\ \text{ for any $B_3$-edge-cut $F'$ of $B_m$}.
\end{equation}

Let $F$ be a minimum $B_3$-edge-cut of $B_n$. Clearly, $B_n-F$ has
exactly two connected components, denoted by $X$ and $Y$,
respectively, and assume $|X|\leq |Y|$ without loss of generality.

Use notations $\delta(X)$ and $\delta(Y)$ to denote the minimum
degrees of $X$ and $Y$, respectively. We assert that $\delta(X)\geq
2$ and $\delta(Y)\geq 2$.

Assume to the contrary that there exists a vertex $x\in X$ such that
$d_X(x)=1$. Then $x$ has at least $3$ neighbors in $Y$ since $B_n$
is $(n-1)$-regular and $n\geq 5$. Let $X'=X\setminus \{x\}$, $Y'=Y
\cup \{x\}$ and $F'$ be the set of edges between $X'$ and $Y'$. Then
$F'$ is also a $B_3$-edge-cut of $B_{n}$ and $|F'| \leq
|F|-3+1=|F|-2$, which contradicts to the minimality of $F$.

To complete our proof, for a fixed $t\in \{1,n\}$ and any $i\in
I_n$, let
 \begin{displaymath}{}
 \begin{array}{ll}
 X_i=X\cap V(B^{t:i}_{n}),&\
 Y_i=Y\cap V(B^{t:i}_{n}),\\
 F_i=F\cap E(B^{t:i}_{n}),&\
 F_{ij}=F\cap E(B^{t:i}_{n},B^{t:j}_{n}),
 \end{array}
\end{displaymath}
where $E(B^{t:i}_{n},B^{t:j}_{n})$ denotes the set of edges between
$B^{t:i}_{n}$ and $B^{t:j}_{n}$ for $i\ne j$. Let
\begin{displaymath}
\begin{array}{l}
 J_X=\{i\in I_n:\ X_i\ne\emptyset\}, J_Y=\{i\in I_n:\ Y_i\not=\emptyset\}\ \ {\rm and} \ \ J_0 =J_X\cap J_Y.
 \end{array}
 \end{displaymath}
Clearly, $|J_X|\geq 1$ and $|J_Y|\geq 1$ since $F$ is a
$B_3$-edge-cut of $B_n$. We choose such $t\in\{1,n\}$ that $|J_X|$
is as large as possible, say $t=1$.

First assume $|J_0|\geq 3$. For $i\in J_0$, we have $\delta(X_i)\geq
1$ and $\delta(Y_i)\geq 1$ since  $\delta(X)\geq 2$ and
$\delta(Y)\geq 2$, and each vertex in $B^{1:i}_{n}$ only has one
neighbor not in $B^{1:i}_{n}$ by Lemma~\ref{lem4.1}. It follows that
$F_i$ is a 1-edge-cut of $B^{1:i}_{n}$ for $i\in J_0$. By
Lemma~\ref{lem4.3}, we have that
   \begin{equation}\label{e4.3}
  |F_i|\geq \lambda^{1}(B_{n-1})=2(n-3) \ \ \text{for}  \  i\in J_0,
  \end{equation}
and so
\[ |F| \geq \sum_{i\in J_0}|F_i| \geq |J_0| 2(n-3)\geq 6(n-3). \]
The inequality (\ref{e4.1}) follows.

Now assume $|J_0|\leq 2$. Let $a=|J_X\setminus J_0|$ and
$b=|J_Y\setminus J_0|$. Then $2a\leq n$ since $|X|\leq |Y|$. Since
any $i\in I_n$ not in $J_X$ is certainly in $J_Y$, $a+b=n-|J_0| \geq
n-2$.

Assume $a\geq 1$ and $j_1\in J_X\setminus J_0,j_2\in J_Y\setminus
J_0$. By Lemma~\ref{lem4.1} there are $(n-2)!$ independent edges
between $B^{1:j_1}_{n}$ and $B^{1:j_2}_{n}$, and so there are
$(ab(n-2)\,!)$ independent edges between $\cup_{j_1\in J_X\setminus
J_0}B^{1:j_1}_{n}$ and $\cup_{j_2\in J_Y\setminus
J_0}B^{1:j_2}_{n}$. It follows that for $n\geq 5$
 $$
\begin{array}{rl}
 |F| &\geq  ab(n-2)\,!
 \geq a(n-2-a)(n-2)\,!\\
 &\geq (n-2-1)(n-2)\,!
 \geq 6(n-3).
\end{array}
$$
The inequality (\ref{e4.1}) follows.

Next assume $a=0$. Then $1\leq |J_X|=|J_0|\leq 2$ and $|J_Y|=n$. We
first show that
\begin{equation}\label{e4.4}
  \text{$Y_i$ contains a $B_3$ as a subgraph for each
 $i\in I_n$}.
 \end{equation}

To this end, for $i\in I_n$, let $K_i$ be the set of the $n$-th
digits of vertices in $X_i$, that is, $K_i=\{p_n:\ ix_2\cdots
x_{n-1}p_n \in X_i,\ p_n\in I_n \setminus \{i\}, x_j\in I_n
\setminus \{i,p_n\}, 2\leq j\leq n-1 \}$. Let $K'_i=I_n\setminus
(K_i\cup\{i\})$. By choice of $t$ and $|J_X| \leq 2$, we have
$|K_i|\leq 2$ and $|K'_i|=(n-1)-|K_i|\geq n-3\geq 2$. Thus, for
fixed $i\in I_n$ and $p_n\in K'_i$, the subgraph induced by
$\{ix_2\cdots x_{n-1}p_n:\ x_j\in I_n \setminus \{i,p_n\}, 2\leq
j\leq n-1\}$ is isomorphic to a $B_{n-2}$ and is contained in $Y_i$,
which implies that $Y_i$ contains a $B_3$ since $n\geq 5$, and so
the assertion (\ref{e4.4}) holds. Thus, for any $i\in J_X$, if $X_i$
contains a $B_3$, then $F_i$ is a $B_3$-edge-cut of $B^{1:i}_n$ by
the assertion (\ref{e4.4}). By the induction hypothesis
(\ref{e4.2}), we have that
 \begin{equation}\label{e4.5}
  \text{$|F_i|\geq 6(n-4)$ if $X_i$ contains a $B_3$ as a subgraph for each
 $i\in J_X$}.
 \end{equation}

We consider two cases according as $|J_X|=1$ and $|J_X|=2$.

{\bf Case 1.}\quad $|J_X|=1$.

Without loss of generality, assume that $J_0=J_X=\{1\}$, and so
$X_1=X\cap V(B^{1:1}_{n})$. By the hypothesis of $F$, $X_1$ contains
a $B_3$, $|X_1|\geq 6$, and each vertex in $X_1$ has one neighbor
not in $B^{1:1}_n$. Combining these facts with the assumption
(\ref{e4.5}), we have that
 $$
 |F|\geq  |F_{1}|+|X_1|
 \geq 6(n-4)+6
 =6(n-3).
$$

{\bf Case 2.}\quad $|J_X|=2$.

Without loss of generality, assume $J_X=\{1,2\}$, and so $X_1=X\cap
V(B^{1:1}_{n})$ and $X_2=X\cap V(B^{1:2}_{n})$. By the hypothesis of
$F$, $X=X_1\cup X_2$ contains a $B_3$. If neither $X_1$ nor $X_2$
contain a $B_3$, then there is a vertex in the $B_3$ whose the first
digit is different from $1$ and $2$, which implies $|J_X|\geq 3$, a
contradiction. Thus, without loss of generality, assume that $X_1$
contains a $B_3$. By the assumptions (\ref{e4.5}) and (\ref{e4.3}),
for $n\geq 6$ we have that
 \begin{equation}\label{e4.6}
 |F|\geq  |F_{1}|+|F_2|
 \geq 6(n-4)+2(n-3)
 \geq 6(n-3).
 \end{equation}

\begin{figure}[h]
\begin{center}
\begin{pspicture}(-2.,-2.2)(2,2.8)
\rput{0}{%
\SpecialCoor\degrees[6]
 \ncline{10}{11}\ncline{11}{12}
 \multido{\i=0+1}{6}{\cnode(.9;\i){.11}{1\i}}
 \ncline{10}{11}\ncline{11}{12}\ncline{12}{13}\ncline{13}{14}\ncline{14}{15}\ncline{15}{10}
 }%
 \cnode(-.45,1.8){.11}{a1}\cnode(.45,1.8){.11}{a2}\ncline{a1}{a2}\ncline{a1}{12}\ncline{a2}{11}
 \cnode(-2.,-.6){.11}{b1}\cnode(-1.6,-1.4){.11}{b2}\ncline{b1}{b2}\ncline{b1}{13}\ncline{b2}{14}
 \cnode(2.,-.6){.11}{c2}\cnode(1.6,-1.4){.11}{c1}\ncline{c1}{c2}\ncline{c1}{15}\ncline{c2}{10}

  \rput(-1.,0.8){\scriptsize$12345$}\rput(-1.4,0.1){\scriptsize$13245$}\rput(-1.,-.7){\scriptsize$13425$}
 \rput(1,0.8){\scriptsize$12435$}\rput(1.4,0.1){\scriptsize$14235$}\rput(1,-.7){\scriptsize$14325$}

  \rput(-.8,2.1){\scriptsize$21345$}\rput(.8,2.1){\scriptsize$21435$}
  \rput(-2.4,-.9){\scriptsize$31245$}\rput(2.4,-.9){\scriptsize$41235$}
  \rput(-2.1,-1.6){\scriptsize$31425$}\rput(2.1,-1.6){\scriptsize$41325$}

 \rput(-1.,2.1){\rnode{a}{}} \rput(1.,2.1){\rnode{b}{}}
 \ncbox[nodesep=8pt,boxsize=.6,linearc=.3]{a}{b}
 \rput(-0.,2.4){\scriptsize$B^{1:2}_5$}
 \rput(-3.,-1.6){\rnode{a}{}} \rput(-1.7,-1.){\rnode{b}{}}
 \ncbox[nodesep=8pt,boxsize=.7,linearc=.3]{a}{b}
 \rput(-2.9,-1.5){\scriptsize$B^{1:3}_5$}
  \rput(3.,-1.6){\rnode{a}{}} \rput(1.7,-1.){\rnode{b}{}}
 \ncbox[nodesep=8pt,boxsize=.7,linearc=.3]{a}{b}
 \rput(2.9,-1.5){\scriptsize$B^{1:4}_5$}

 \end{pspicture}
\caption{\label{f5}\footnotesize {The neighbors of a $B_3$ in a
$B_5$} }
\end{center}
\end{figure}
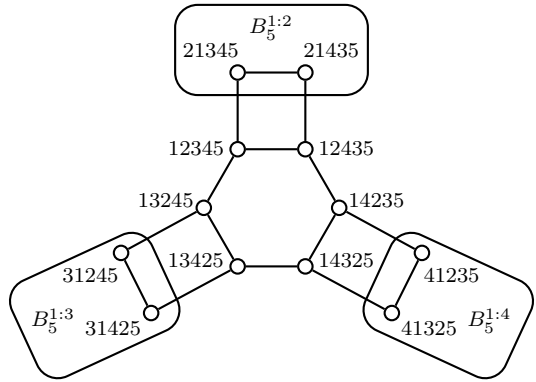

Now assume $n=5$. Then the vertices of $B_3$ in $X_1$ must have
forms $1x_2x_3x_4p_5$ or $1p_2x_3x_4x_5$, where $p_2$ and $p_5$ are
fixed. Since $|K_i|\leq 2$ for each $i\in I_n$, the vertices of
$B_3$ in $X_1$ have the former form. Without loss of generality, say
$p_5=5$. Then the vertex-set of $B_3$ and the neighbors of $X_1$ not
in $B^{1:1}_5$ are shown in Fig.~\ref{f5}, where four neighbors are
in $B^{1:3}_5\cup B^{1:4}_5$. This shows that no matter how
$x_2x_3x_4p_5$ is chosen, we always have $\sum^5_{j=3}F_{1j}\geq
2\times 2=4$. It follows that
$$
 |F| \geq  |F_{1}|+|F_2|+ \sum^n_{j=3}F_{1j}
 \geq 6(n-4)+2(n-3)+4
 >6(n-3).
$$

The inequality (\ref{e4.1}) follows. By induction principles, the
theorem follows.
\end{pf}

\section{Conclusions}

In this paper, we investigate the $h$-embedded connectivity
$\zeta_h$ and $h$-embedded edge-connectivity $\eta_h$ in the
hypercube $Q_n$, the star graph $S_n$ and the bubble-sort graph
$B_n$. We determine that $\zeta_h(Q_n)=2^h(n-h)$ for $h\leq n-2$,
$\eta_h(Q_n)=2^h(n-h)$ for $h\leq n-1$,
$\zeta_{h}(S_{n})=\eta_{h}(S_{n})=h!(n-h)$ for $1 \leq h \leq n-1$,
and $\eta_3(B_n)=6(n-3)$ for $n\geq 4$. These results can provide
more accurate measurements for fault tolerance of the system when
the graphs are used to model the topological structure of
large-scale parallel processing systems, and correct some wrong
conclusions and proofs in~\cite{yw12,yw14}. The value of
$\zeta_h(B_n)$ for $h\geq 3$ and the value of $\eta_h(B_n)$ for
$h\geq 4$ deserve further research.


\vskip6pt

\begin{thebibliography}{10}

\bibitem{ak89}
S. B. Akers, B. Krishnamurthy, A group theoretic model for symmetric
interconnection networks. IEEE Transactions on Computers, 38 (4)
(1989), 555-566.


\bibitem{h83}
F. Harary, Conditional connectivity. Networks, 13 (1983), 347-357.


\bibitem{lhm94}
S. Latifi, M. Hegde, M. Naraghi-Pour, Conditional connectivity
measures for large multiprocessor systems. IEEE Transactions on
Computers, 43 (2) (1994), 218-222.

\bibitem{lx14}
X.-J. Li, J.-M. Xu, Generalized measures for fault tolerance of star
networks, Networks, 63 (3) (2014), 225-230.


\bibitem{oc93}
A. D. Oh, H. Choi, Generalized measures of fault tolerance in
$n$-cube networks. IEEE Transactions on Parallel and Distributed
Systems, 4 (1993), 702-703.

\bibitem{sls12}%
W. Shi, F. Luo, P. K. Srimani, A new hierarchical structure of star
graphs and applications.
Lecture Notes in Computer Science, 7154, (2012), pp. 267-268.

\bibitem{wg98}
J. Wu, G. Guo, Fault tolerance measures for $m$-ary $n$-dimensional
hypercubes based on forbidden faulty sets. IEEE Transactions on
Computers, 47 (1998), 888-893.

\bibitem{x00c}
J.-M. Xu, On conditional edge-connectivity of graphs. Acta
Mathematicae Applicatae Sinica, 16 (4) (2000), 414-419.

\bibitem{x00d}
J.-M. Xu, Restricted edge-connectivity of vertex-transitive graphs,
Chinese Journal of Contemporary Mathematics, 21 (4) (2000), 369-374.

\bibitem{x01}
J.-M. Xu, Topological Structure and Analysis of Interconnection
Networks. Kluwer Academic Publishers, Dordrecht/Boston/London, 2001.


\bibitem{ylm10}
W.-H. Yang, H.-Z Li, J.-X Meng, Conditional connectivity of Cayley
graphs generated by transposition trees. Information Processing
Letters, 110 (2010), 1027-1030.

\bibitem{yw12}
Y. Yang,  S. Wang, Conditional connectivity of star graph networks
under embedded restriction. Information Sciences,  199 (15) (2012),
187-192.

\bibitem{yw14}
Y. Yang, S. Wang, J. Li, Conditional connectivity of recursive
interconnection networks respect to embedding restriction.
Information Sciences, 279 (2014), 273-279.

\end{thebibliography}
\end{document}